\documentclass[12pt]{article}
\usepackage{amsmath}
\usepackage{amssymb}
\usepackage{amscd}
\newtheorem{theoreme}{Th\'eor\`eme}
\newtheorem{proposition}{Proposition}
\newtheorem{corollaire}{Corollaire}
\newtheorem{lemme}{Lemme}

\newtheorem{remarque}{Remarque}
\newtheorem{question}{Question}
\newenvironment{preuve}{\begin{trivlist} \item[]{\it Preuve---}}
{\par\hfill $\square$\end{trivlist}}
\renewcommand{\P}{\mathbb{P}}
\newcommand{\C}{\mathbb{C}}
\newcommand{\N}{\mathbb{N}}
\newcommand{\R}{\mathbb{R}}
\newcommand{\T}{{\rm T}}
\newcommand{\id}{{\rm id}}
\newcommand{\D}{{\cal D}}
\newcommand{\F}{{\cal F}}

\title{Ensembles d'unicit\'e pour les polyn\^omes
\footnote{{\bf Classification math\'ematique:} 30D05, 58F23.
\ \ \ \ \ \ \ \ \ \ \ \ \ \ \ \ \ \ \ \ \ \ \ \ \ \ \
\break
{\bf Mots cl\'es:} ensemble d'unicit\'e, 
ensemble de Julia, mesure invariante.}
}
\author{Tien-Cuong Dinh
\footnote{Math\'ematique-B\^atiment 425,
Universit\'e Paris-Sud, 91405 ORSAY Cedex (France).\break
E-mails: TienCuong.Dinh@math.u-psud.fr.
}}
\begin{document}
\maketitle
\begin{abstract} Let $E\subset \C$ 
be a compact set of positive logarithmic capacity. Let us suppose that
  for every polynomial $P\not=\id$ we have $P^{-1}(E)\not=E$. Then for 
all no constant polynomials $f$ and $g$ such that
  $f^{-1}(E)=g^{-1}(E)$ we have $f=g$. 
\end{abstract}
\section{Introduction}
On note $\P^1=\C\cup\{\infty\}$ la droite projective complexe. 
Un compact $E$ de $\C$
est appel\'e {\it ensemble d'unicit\'e} si pour tous
polyn\^omes non constants
$f$ et $g$ v\'erifiant $f^{-1}(E)=g^{-1}(E)$ on a
$f=g$. S'il existe un polyn\^ome $P\not=\id$ tel que
$P^{-1}(E)=E$, alors $E$ n'est pas un ensemble d'unicit\'e.
\begin{question} Supposons que $P^{-1}(E)\not=E$ pour tout
polyn\^ome $P\not=\id$. $E$ est-t-il un ensemble
d'unicit\'e?   
\end{question}
Les ensembles d'unicit\'e pour les polyn\^omes de m\^eme degr\'e
sont d\'etermin\'es par Ostrovskii, Pakovitch et Zaidenberg
\cite{OPZ, Pakovitch}. Les ensembles
d'unicit\'e pour les fonctions enti\`eres ou m\'eromorphes avec
un nombre minimal d'\'el\'ements sont \'etudi\'es par Nevanlinna
et \'egalement par plusieurs
autres auteurs ({\it voir} par exemple \cite{GrossYang,
Shiffman, KhoaiKhue}).
\begin{question}
Soit $\mu$ une mesure de probabilit\'e \`a support compact dans
$\C$.
Pour quels polyn\^omes $f$ et $g$ de degr\'es
$d\geq 1$ et $d'\geq 1$ on a $d^{-1}f^*(\mu)={d'}^{-1}g^*(\mu)$?
\end{question}
Les deux questions pos\'ees ci-dessus
sont \'etroitement li\'ees.
En effet, si la
capacit\'e logarithmique $E$ est positive et si
$f^{-1}(E)=g^{-1}(E)$, alors $d^{-1}f^*(\mu)={d'}^{-1}g^*(\mu)$ pour
la mesure d'\'equilibre $\mu$ de $E$. R\'eciproquement,
si $E$ est le
support de $\mu$ et si $d^{-1}f^*(\mu)={d'}^{-1}g^*(\mu)$,
on a $f^{-1}(E)=g^{-1}(E)$. De plus, on peut
trouver un compact $E_0$ de capacit\'e logarithmique positive
tel que $f^{-1}(E_0)=g^{-1}(E_0)$ ({\it voir} le paragraphe 3). 
Lorsque $f\not=g$,
$E$ et $E_0$ ne sont donc pas ensembles d'unicit\'e.
Les deux cas particuliers de ces
probl\`emes sont le probl\`eme de
d\'etermination des fonctions ayant le m\^eme ensemble de Julia ou
la m\^ eme mesure totalement invariante et le probl\`eme de
d\'etermination des fonctions permutables \cite{Fatou, Julia,
BakerEremenko} ({\it voir} \'egalement
\cite{Ritt, Eremenko, LevinPrzytycki, Dinh1, Dinh2}).
Notre r\'esultat principal
est le th\'eor\`eme suivant:
\begin{theoreme} Soient $\mu$ une mesure de probabilit\'e \`a
support compact dans $\C$ et $f$,
$g$ deux polyn\^omes de degr\'es $d\geq 1$ et $d'\geq 1$. Soit
$m$ le plus grand diviseur commun de $d$ et $d'$.
Supposons que $d^{-1}f^*(\mu)={d'}^{-1}g^*(\mu)$. Alors il existe un
polyn\^ome $Q$ de degr\'e $m$ et des polyn\^omes $f_0$, $g_0$
tels que $f=f_0\circ Q$, $g=g_0\circ Q$ et tels que l'une des
conditions suivantes soit vraie:
\begin{enumerate}
\item $f_0=\id$ ou $g_0=\id$.
\item $d>m$, $d'>m$ et
pour une certaine coordonn\'ee $z$ de $\C$ on a $f_0(z)=z^{d/m}$,
$g_0(z)= az^{d'/m}$ o\`u $a\not = 0$ est une constante.
\item $d>m$, $d'>m$ et
pour une certaine coordonn\'ee $z$ de $\C$ on a $f_0=\pm \T_{d/m}$,
$g_0=\pm \T_{d'/m}$ o\`u $\T_k$ est le polyn\^ome de Tchebychev
de degr\'e $k$. 
\end{enumerate}   
\end{theoreme}
\begin{corollaire} Soit $E$ un compact de capacit\'e
logarithmique positive de
$\C$. Alors $E$ est un ensemble d'unicit\'e si et seulement si
pour tout polyn\^ome $P\not=\id$ on a $P^{-1}(E)\not=E$.
\end{corollaire}
Soit $P$ un polyn\^ome de degr\'e au moins deux. {\it L'ensemble
de Julia rempli} $K_P$
de $P$ est l'ensemble des points d'orbite
born\'ee, i.e. les points $z$ tels que les suites
$\{P^n(z)\}_{n\in \N}$ soient born\'ees. Ici on note
$P^n:=P\circ\cdots \circ P$ le $n$-i\`eme
it\'er\'e de $P$. {\it L'ensemble de Julia} $J_P$ est
le bord topologique de l'ensemble $K_P$. Alors $K_P$
est le plus grand compact {\it totalement invariant} par $P$, i.e.
$P^{-1}(K_P)=K_P$. L'ensemble $J_P$ est le
plus petit compact totalement invariant par $P$ qui contient
plus qu'un \'el\'ement ({\it voir} par exemple \cite[4.2.2]{Beardon}).
\begin{corollaire} Soit $E$ un compact de capacit\'e
logarithmique positive de
$\C$. Supposons que $E$ n'est pas un ensemble d'unicit\'e et que
$E$ n'est pas invariant par aucune rotation de $\C$. Alors il
existe un polyn\^ome $P$ de degr\'e au moins deux tel que
$P^{-1}(E)=E$ et $J_P\subset E\subset K_P$.
\end{corollaire}
Si dans un ouvert $U$ l'intersection  $J_P\cap U$ est un ensemble
non vide, inclus dans une courbe r\'eelle lisse, alors $P(z)=z^d$
ou $P(z)=\pm\T_d(z)$ pour une certaine coordonn\'ee $z$ de $\C$
\cite[p.127]{Steinmetz}. La
notation $\T_d$ signifie {\it le polyn\^ome de Tchebychev} de degr\'e
$d$ d\'efini par $\T_d(\cos t):=\cos dt$. Si $P(z)=z^d$, $K_P$
est le disque unit\'e, $J_P$ est le cercle unit\'e; si
$P(z)=\pm \T_d(z)$, $K_p$ et $J_p$ sont \'egaux au segment $[-1,1]$.
On en d\'eduit facilement que,
toute courbe r\'eelle lisse par morceaux,
qui n'est invariante par aucune rotation, est un
ensemble d'unicit\'e.
\par
L'ensemble $J_P$ est le support d'une mesure de probabilit\'e
$\mu_P$ qui est {\it totalement
invariante} par $P$, i.e. $(\deg P)^{-1}P^*(\mu_P)=\mu_P$. C'est
la seule mesure de probabilit\'e \`a support compact
qui est
totalement invariante par $P$ sauf dans le cas o\`u $P(z)=z^d$
pour une certaine
coordonn\'ee $z$ de $\C$. Dans ce cas
exceptionnel, toute mesure totalement invariante est une
combinaison lin\'eaire de la mesure de Lebesgue sur le cercle
unit\'e et de la masse de Dirac en $0$. Dans le
th\'eor\`eme 1, si $g_0=\id$, $\deg f_0\geq 2$ et si $f_0$ n'est
pas conjugu\'e \`a $z^{d/m}$ on a
$\mu=\mu_{f_0}$. Dans la deuxi\`eme condition du th\'eor\`eme 1, $\mu$
est une combinaison lin\'eaire de la masse de Dirac en $0$ et de la
mesure de Lebesgue sur le cercle $\{z\in\C: |z|=|a|^{d/(d-d')}\}$. 
Dans la troisi\`eme
condition de ce th\'eor\`eme, $\mu$ est la mesure totalement
invariante des polyn\^omes de Tchebychev.
\par
Dans la preuve du th\'eor\`eme 1, on se ram\`ene \`a des syst\`emes
dynamiques holomorphes en une et en plusieurs variables. D'abord,
on peut choisir une fonction
subharmonique  $\varphi_0$ telle que
$i\partial\overline\partial\varphi_0=\mu$ et $d^{-1}\varphi_0\circ
f={d'}^{-1}\varphi_0\circ g=:\varphi_{-1}$.
Notons $\delta_f(z):=\exp(2i\pi/d)
z+a_0+a_1z^{-1}+\cdots$ le germe
d'application holomorphe d\'efini au
voisinage de $\infty$ qui pr\'eserve les fibres de $f$. Notons
$\delta_g$ le germe analogue pour $g$. Alors $\delta_f$ et
$\delta_g$ pr\'eservent les lignes de niveau de $\varphi_{-1}$. Dans une
coordonn\'ee locale convenable, les
lignes de niveau de $\varphi_{-1}$ au voisinage de $\infty$ sont
les cercles de centre $\infty$. Par cons\'equent, pour cette
coordonn\'ee locale, $\delta_f$ et $\delta_g$ sont des rotations.
D'o\`u 
$\delta_f\circ\delta_g=\delta_g\circ\delta_f$ et
$\delta_f^{d/m}=\delta_g^{d'/m}=:\delta$. Notons ${\cal F}_f(z)$
{\it la fibre} de
$f$ qui contient $z$, i.e. ${\cal F}_f(z)=f^{-1}\circ f(z)$. 
Le polyn\^ome $Q$ sera d\'efini
comme un polyn\^ome dont toute fibre g\'en\'erique
$\F_Q(z)$ est \'egale \`a
l'intersection $\F_f(z)\cap \F_g(z)$.
Au voisinage de $\infty$, la fibre 
$\F_Q(z)$ est l'orbite de $z$ par $\delta$.
Ceci entra\^{\i}ne que $Q$ est
bien d\'efini et qu'il est de degr\'e $m$ (proposition 2). 
Si $m=d$ ou $m=d'$, on peut choisir
$P=f$ ou $P=g$; la condition 1 du th\'eor\`eme 1 est alors vraie.
Pour la suite de la preuve,
on peut supposer que $d>1$ et $d'>1$ sont
premiers entre eux.
\par
Notons $\Phi$ un polyn\^ome de degr\'e $dd'$ v\'erifiant
$\F_\Phi(z)=g^{-1}\circ g(\F_f(z))$ pour tout $z\in \C$.
Au voisinage de $\infty$,
$\F_\Phi(z)$ est l'ensemble des points 
$\delta_g^n\circ\delta_f^m(z)$ car
$\delta_f$ et $\delta_g$ commutent. Ceci entra\^{\i}ne que $\Phi$
existe et que
$\F_\Phi(z)=f^{-1}\circ
f(\F_g(z))$ (proposition 2). Alors on peut
d\'ecomposer $\Phi=f_1\circ g=g_1\circ f$ o\`u $f_1$ (resp. $g_1$) est
un polyn\^ome de degr\'e $d$ (resp. $d'$). En suite, on peut montrer 
que ${\cal C}_{f_1}=g({\cal C}_f)$ et ${\cal C}_{g_1}=f({\cal
C}_g)$ o\`u ${\cal C}$ signifie l'ensemble critique. Ceci, sous
certaines conditions pos\'ees sur les coefficients dominants et
sur les valeurs de polyn\^omes en $0$, permet de construire un
endomorphisme polynomial $\D_{d,d'}$ de l'ensemble $\Sigma(d,d')$
des couples $(f,g)$. L'ensemble des couples $(f,g)$ qui v\'erifient
les hypoth\`eses du th\'eor\`eme 1, est un sous-ensemble
alg\'ebrique ${\cal N}(d,d')$ invariant par $\D_{d,d'}$.
L'ensemble des couples $(f,g)$ qui
v\'erifient la condition 2 ou 3 du
th\'eor\`eme 1, d\'ecrit deux courbes ${\cal C}_1(d,d')$ et
${\cal C}_2(d,d')$. Il reste \`a prouver que ${\cal
N}(d,d')={\cal C}_1(d,d')\cup{\cal C}_2(d,d')$. On 
montrera que les points p\'eriodiques de ${\cal N}(d,d')$
appartiennent \`a ${\cal C}_1(d,d')\cup{\cal C}_2(d,d')$. Ceci
est d\^u \`a la solution d'une \'equation bien connue: $f\circ
g=g\circ f$ (en particulier, pour les points fixes, on obtient
directement $\Phi=f\circ g=g\circ f$).
Finalement, l'invariance de ${\cal
N}(d,d')$ par $\D_{d,d'}$ implique que ${\cal
N}(d,d')={\cal C}_1(d,d')\cup{\cal C}_2(d,d')$.
\par
{\it Remerciement.}--- Je tiens \`a remercier Charles Favre et
Nessim Sibony pour leurs aides pendant la pr\'eparation de cet
article. 
\section{Factorisation de polyn\^omes}
Soit $P$ un polyn\^ome de degr\'e $d\geq 1$.
Il existe une fonction holomorphe
unique $B(z)=z+a_0+a_1z^{-1}+a_2z^{-2}+\cdots$ d\'efinie au
voisinage de $\infty$ telle que $B\circ P\circ B^{-1}=\alpha
z^d$ o\`u $\alpha$ est le coefficient
dominant de $P$ ({\it voir} par exemple \cite[6.10.1]{Beardon}). On
d\'efinit la fonction $\delta_P$ par la formule
$\delta_P(z):= B^{-1}(\theta_d B(z))$ o\`u
$\theta_d:=\exp(2\pi i/d)$. Cette fonction est
d\'efinie au voisinage de $\infty$, permute les \'el\'ements de
chaque fibre de $P$ et v\'erifie $P\circ \delta_P=P$, 
$\delta_P^d=\id$. Soit $\Phi$
une application d\'efinie dans un voisinage suffisamment petit
de $\infty$ \`a l'image
dans un espace $X$. Alors $\Phi$ s'\'ecrit sous la forme
$G\circ P$ si et seulement si $\delta_P$ pr\'eserve les
fibres de $\Phi$, i.e. $\Phi\circ\delta_P=\Phi$.
\begin{proposition} Soient $X$ un espace m\'etrique, $a$ un
point de $X$ et $\Phi$ une application d\'efinie
dans un voisinage suffisamment petit
de $\infty$ \`a l'image dans $X\setminus\{a\}$
v\'erifiant
$\lim_{z\rightarrow \infty} \Phi(z)=a$. Soient $P_1$
et $P_2$ deux polyn\^omes de degr\'es $d_1\geq 1$ et $d_2\geq 1$
tels que $\Phi$ s'\'ecrive sous les formes
$\Phi=G_1\circ P_1=G_2\circ P_2$. Soit $m$ le plus grand diviseur
commun de $d_1$ et $d_2$. Alors
$\delta_{P_1}\circ\delta_{P_2}=
\delta_{P_2}\circ\delta_{P_1}$ et $\delta_{P_1}^{d_1/m} =
\delta_{P_2}^{d_2/m}$. En particulier, si $d_1=d_2$ il existe un
polyn\^ome lin\'eaire $\sigma$ tel que $P_1=\sigma\circ P_2$.
\end{proposition}
\begin{preuve} 
On montre
que $\delta_{P_1}\circ\delta_{P_2}=
\delta_{P_2}\circ\delta_{P_1}$.
Supposons que
$\delta:=\delta_{P_1}^{-1}\circ \delta_{P_2} \circ
\delta_{P_1}\circ \delta_{P_2}^{-1}\not =\id$.
Remarquons que $\delta$ s'\'ecrit sous la forme
$\delta(z)=z+a_0+a_1z^{-1}+a_2z^{-2}+\cdots$. La dynamique d'une
telle application est bien connue ({\it voir} par exemple
\cite[6.5]{Beardon}). Il existe
un point $z$ tel que $\delta^n(z)$ tende vers $\infty$ quand
$n\rightarrow +\infty$. Par cons\'equent, $\Phi(\delta^n(z))$ tend
vers $a$. C'est une contradiction car $\Phi\circ\delta=\Phi$. De
m\^eme mani\`ere, on montre que $\delta_{P_1}^{d_1/m}\circ
\delta_{P_2}^{-d_2/m}= \id$ et
$\delta_{P_1}^{d_1/m}=\delta_{P_2}^{d_2/m}$.
\par
Si $d_1=d_2$, on a $\delta_{P_1}=\delta_{P_2}$. Par cons\'equent,
$\F_{P_1}(z)=\F_{P_2}(z)$ au voisinage de $\infty$. Par analyticit\'e,
ceci est vrai pour tout $z\in\C$. Alors on peut d\'efinir la fonction
$\sigma$ holomorphe dans $\C$ par
$\sigma:=P_1\circ P_2^{-1}$. Il est clair que $\sigma$ est
bijective. Par cons\'equent, $\sigma$ est lin\'eaire et
$P_1=\sigma\circ P_2$.
\end{preuve}
Pour tout polyn\^ome $P$, on note ${\cal C}_P$ l'ensemble critique de
$P$. Un point 
de ${\cal C}_P$ sera compt\'e $k$ fois s'il est de multiplicit\'e $k$.
\begin{proposition} Soient $P_1$ et $P_2$ deux polyn\^omes 
de degr\'es $d_1\geq 1$ et $d_2\geq 1$. 
Soit $m$ le plus grand diviseur commun de $d_1$ et $d_2$. 
\begin{enumerate}
\item Si $\Phi$ est un polyn\^ome v\'erifiant
$\Phi\circ \delta_{P_1}=\Phi$ au voisinage de $\infty$, 
alors il existe un polyn\^ome $R$ tel que $\Phi=R\circ P_1$.
\item Si $\delta_{P_1}^{d_1/m}=\delta_{P_2}^{d_2/m}$ 
alors il existe un polyn\^ome $Q$ de 
degr\'e $m$ et des polyn\^omes $R_1$, $R_2$ tels que 
$P_1= R_1\circ Q$ et $P_2=R_2\circ Q$. 
En particulier, si $d'=m$ il existe un polyn\^ome $R$ 
tel que $P_1=R\circ P_2$.
\item Si $m=1$ et si $\delta_{P_1}\circ 
\delta_{P_2}=\delta_{P_2}\circ\delta_{P_1}$, 
il existe un polyn\^ome $\Phi$ de degr\'e $d_1d_2$ 
et des polyn\^omes $P_1^*$, $P_2^*$ tels que 
$\Phi=P_1^*\circ P_2=P_2^*\circ P_1$. De plus, 
on a ${\cal C}_{P_1^*}=P_2({\cal C}_{P_1})$
et ${\cal C}_{P_2^*}=P_1({\cal C}_{P_2})$.
\end{enumerate}
\end{proposition}
\begin{preuve}
{\bf 1.} Le fait 
que $\Phi\circ \delta_{P_1}$ implique que $\F_{P_1}(z)\subset 
\F_{\Phi}(z)$ pour $z$ dans un voisinage de $\infty$. 
Par analyticit\'e, ceci est vrai pour tout $z$. 
Par cons\'equent, on peut d\'efinir la fonction $R$ holomorphe 
dans $\C$ par $R(z):=\Phi\circ P^{-1}(z)$. On a $\Phi=R\circ 
P_1$. Comme $\Phi$ et $P$ 
sont des polyn\^omes, $\lim_{z\rightarrow \infty} R(z)=\infty$. 
Donc $R$ est un polyn\^ome. 
\par 
{\bf 2.} On note $w:=(w_1,w_2)$ les coordonn\'ees de $\C^2$. 
Soient $\Pi :\C\longrightarrow \C^2$ d\'efini par 
$\Pi(z):=(P_1(z),P_2(z))$ et ${\cal C}:=\Pi(\C)$. Alors  
${\cal C}$ est une courbe alg\'ebrique de $\C^2$. De plus, elle 
est parabolique car elle est l'image de $\C$ par une application 
holomorphe. D'o\`u ${\cal C}$ est un $\C$ ou un $\C^*$ immerg\'e 
dans $\C^2$. La courbe compactifi\'ee $\overline{\cal 
C}\subset\P^2$ est donc rationnelle. Comme $P(z)/Q(z)$ a une 
limite finie ou infinie quand $z\rightarrow \infty$, 
$\overline{\cal C}$ coupe la droite infinie $L$ en un seul point 
$a$. De plus, au voisinage de $a$, $\overline{\cal C}$ est 
irr\'eductible car c'est l'image d'un voisinage de $\infty$ par 
l'application $\Pi$. On en d\'eduit que ${\cal C}=\overline{\cal 
C}\setminus \{a\}$ est un $\C$ immerg\'e. Soit $\varphi: 
\C\longrightarrow {\cal C}$ une application holomorphe, injective en 
dehors d'un nombre fini de points. On pose 
$Q:=\varphi^{-1}\circ\Pi$. Alors $Q$ est une fonction holomorphe 
d\'efinie en dehors d'un nombre fini de points au voisinage 
desquels elle est born\'ee. Par cons\'equent, $Q$ se prolonge en 
une fonction holomorphe sur $\C$. On v\'erifie facilement que 
$\lim_{z\rightarrow \infty} Q(z)=\infty$. Donc $Q$ est un 
polyn\^ome. On pose 
$\delta:=\delta_{P_1}^{d_1/m}=\delta_{P_2}^{d_2/m}$. 
\par 
Au voisinage de $\infty$, on a 
$$Q\circ\delta=\varphi^{-1}\circ 
(P_1\circ\delta, P_2\circ\delta)=\varphi^{-1}\circ 
(P_1,P_2) =Q.$$ 
\par 
Soient $z_1$ et $z_2$ suffisamment 
proches de $\infty$ tels que $Q(z_1)=Q(z_2)$. Alors 
$P_1(z_1)=P_1(z_1)$ et $P_2(z_1)=P_2(z_2)$. Il existe donc les 
entiers $0\leq n_1\leq d_1-1$ et $0\leq n_2\leq d_2-1$ tels que 
$z_1=\delta_{P_1}^{n_1}(z_2)=\delta_{P_2}^{n_2}(z_2)$. De plus, 
on sait que 
$$\lim_{z\rightarrow\infty} \frac{\delta_{P_j}^{n_j}(z)}{z} 
=\exp(2n_j\pi i/d_j).$$ 
Par cons\'equent, l'\'egalit\'e 
$\delta_{P_1}^{n_1}(z_2)=\delta_{P_2}^{n_2}(z_2)$ pour $z_2$ 
suffisamment proche de $\infty$ implique que $n_jm$ est  
divisible par $d_j$ pour $j=1$ ou $2$. 
On a alors $z_1=\delta^{n_1m/d_1}(z_2)$. 
\par 
Les deux arguments ci-dessus montrent que $\delta_Q=\delta$. Par 
cons\'equent, $\deg Q=m$. D'apr\`es la premi\`ere partie, il 
existe des polyn\^omes $R_1$ et $R_2$ tels que $P_1=R_1\circ Q$ 
et $P_2=R_2\circ Q$. 
\par 
On remarque que $\F_Q(z)=\F_{P_1}(z)\cap\F_{P_2}(z)$ pour un 
$z$ g\'en\'erique car ceci est vrai au voisinage de $\infty$. On 
peut \'egalement prouver cette partie par la m\^eme m\'ethode que l'on
utilisera dans la troisi\`eme partie. 
\par 
{\bf 3.} Comme $m=1$, il existe des entiers relatifs 
$n_1$ et $n_2$ tels que 
$n_1d_2+n_2d_1=1$. Posons 
$\delta:=\delta_{P_1}^{n_1}\circ\delta_{P_2}^{n_2}$. Alors $\delta$ 
s'\'ecrit sous la forme $\delta(z)=\exp(2\pi 
i/dd')z+a_0+a_1z^{-1}+\cdots$.   
Notons $\F(z)=P_2^{-1}\circ P_2(\F_{P_1}(z))$. Un point 
de $\F(z)$ est compt\'e $k$ fois s'il est de multiplicit\'e $k$. 
Alors $\F(z)$ est de cardinal $dd'$ et 
au voisinage de $\infty$, $\F(z)$ est l'orbite de $z$ par 
$\delta$ car $\delta_{P_1}\circ \delta_{P_2}=\delta_{P_2}\circ 
\delta_{P_1}$. On en d\'eduit que $\F(z)=P_1^{-1}\circ 
P_1(\F_{P_2}(z))$ au voisinage de $\infty$. Par analyticit\'e, 
ceci est vrai pour tout $z$. 
\par 
Pour tout $1\leq n\leq dd'-1$, 
on note $S_n(z)$ la somme sym\'etrique des termes du type 
$z_1\ldots z_n$ avec $\{z_1,\ldots,z_n\}\subset\F(z)$. Alors 
$S_n(z)$ est une fonction holomorphe sur $\C$. Il est clair que 
$|\delta_{P_1}^i\circ\delta_{P_2}^j(z)/z|$ tend vers $1$ quand 
$z\rightarrow \infty$ pour tous 
$i$ et $j$. Par cons\'equent, $S_n(z)=\mbox{O}(|z|^n)$ quand 
$z\rightarrow\infty$. On en d\'eduit que $\deg S_n\leq n\leq dd'-1$. 
Comme au voisinage de $\infty$, 
$\F(z)$ est l'orbite de $z$ par $\delta$, on a $\F(z)=\F(z_1)$ 
pour tout $z_1\in\F(z)$. Par cons\'equent, $S_n$ est constant sur 
$\F(z)$. Le fait que $\F(z)$ est de cardinal $dd'>\deg S_n$ 
implique que $S_n$ est un polyn\^ome constant. 
Soit $\Phi$ le polyn\^ome de 
degr\'e $dd'$ d\'efini par: 
$$\Phi(z):=z^{dd'}-S_1 z^{dd'-1}+\cdots+(-1)^{dd'-1}S_{dd'-1}z.$$ 
Alors au voisinage de $\infty$, $\F_\Phi(z)=\F(z)$. Par 
analyticit\'e, ceci est vrai pour tout $z$. Par cons\'equent, 
$\Phi\circ\delta_{P_1}=\Phi$ et $\Phi\circ\delta_{P_2}=\Phi$. 
D'apr\`es la premi\`ere partie, il existe des polyn\^omes $P_1^*$ 
et $P_2^*$ tels que $\Phi=P_1^*\circ P_2=P_2^*\circ P_1$. 
\par 
On a ${\cal C}_\Phi={\cal C}_{P_2}\cup P_2^{-1}({\cal C}_{P_1^*})$. 
Ici la notation ${\cal C}$ signifie l'ensemble critique et un 
point critique sera compt\'e $k$ fois s'il est de multiplicit\'e 
$k$. Remarquons qu'un point critique $z$ 
de $\Phi$ est de multiplicit\'e $k$ si et seulement si 
$z$ est un point 
de multiplicit\'e $k+1$ de $\F(z)$. Comme 
$\F(z)=P_2^{-1}\circ P_2({\cal F}_{P_1}(z))$, on a ${\cal 
C}_\Phi={\cal C}_{P_2}\cup P_2^{-1}\circ P_2({\cal C}_{P_1})$. 
Alors $P_2^{-1}({\cal C}_{P_1^*})=P_2^{-1}\circ P_2 
({\cal C}_{P_1})$  car ${\cal C}_\Phi={\cal C}_{P_2} 
\cup P_2^{-1}({\cal C}_{P_1^*})$. 
D'o\`u ${\cal C}_{P_1^*}=P_2({\cal C}_{P_1})$. 
De m\^eme, on a ${\cal C}_{P_2^*}=P_1({\cal C}_{P_2})$. 
\end{preuve}
Soient $\mu$ une mesure de probabilit\'e \`a support compact
et $f$, $g$ deux
polyn\^ omes de degr\'es $d$ et $d'$
v\'erifiant $d^{-1}f^*(\mu)={d'}^{-1}g^*(\mu)$. Soit
$\varphi$ une fonction subharmonique v\'erifiant
$i\partial\overline \partial \varphi=\mu$. 
Cette fonction est
harmonique sur la composante non born\'ee de
$\C\setminus\mbox{supp}(\mu)$; elle
est unique \`a une constante
pr\`es. De plus, $\varphi-\ln|z|$ est harmonique
et born\'ee au voisinage de $\infty$.
Posons $\psi:=d^{-1}\varphi\circ f$. Alors
$i\partial\overline \partial\psi=d^{-1}f^*(\mu)$. Comme
$d^{-1}f^*(\mu)={d'}^{-1}g^*(\mu)$, on a
${d'}^{-1}\varphi\circ g=\psi+c$ o\`u $c$ est une
constante. Lorsque $d\not=d'$, quitte \`a remplacer $\varphi$
par $\varphi-dd'c/(d-d')$, on peut supposer que $c=0$. Dans tous
les cas, on peut appliquer la proposition 1 pour la
fonction $\psi$. On obtient
$\delta_f\circ\delta_g=\delta_g\circ\delta_f$ et 
$\delta_f^{d/m}=\delta_g^{d'/m}$ o\` u $m$ est le plus grand diviseur
commun de $d$ et $d'$. D'apr\`es la proposition 2, on a:
\begin{corollaire}
Il existe un polyn\^ome $Q$ de degr\'e $m$ et des
polyn\^omes $f_0$, $g_0$ tels que
$f=f_0\circ Q$,
$g=g_0\circ Q$ et $(d/m)^{-1}f_0^*(\mu)=(d'/m)^{-1}g_0^*(\mu)$.
En particulier, si $d'$ divise $d$, il
existe un polyn\^ome $P$ de degr\'e $d/d'$
tel que $f=P\circ g$ et $d^{-1}d'P^*(\mu)=\mu$.  
\end{corollaire}
\section{Endomorphisme polynomial $\D_{d,d'}$}
Soient $\mu_0$ une mesure de probabilit\'e \`a support compact de
$\C$ et
$f_0$, $g_0$ deux
polyn\^omes de degr\'es $d>1$, $d'>1$ v\'erifiant
$d^{-1}f_0^*(\mu_0)={d'}^{-1}g_0^*(\mu_0)$ . On suppose
que $d$ et $d'$ sont premiers entre eux ({\it voir}
 le corollaire 3) et que $d>d'$. 
Soient $\beta\not=0$ et $\alpha\not=0$
les coefficients dominants de $f_0$ et $g_0$.
On choisit un point
$a$ tel que $f_0(a)=g_0(a)$. Soit $b:=f_0(a)=g_0(a)$.
Quitte \`a remplacer $f_0$ par $\sigma_1\circ f_0\circ \sigma_2$,
$g_0$ par $\sigma_1\circ g_0\circ \sigma_2$ et $\mu_0$ par
$(\sigma_1)_*(\mu_0)$ on
peut supposer que $a=b=0$ et $\beta=1$ o\`u
$\sigma_2(z):=Az+a$,
$\sigma_1(z):=A^{-d}\beta^{-1}(z-b)$ et $A\in\C^*$.
Soit $\varphi$ la fonction subharmonique v\'erifiant
$i\partial\overline\partial \varphi=\mu_0$ et
$d^{-1}\varphi\circ f_0={d'}^{-1}\varphi\circ g_0$. Posons
$\varphi_0(z):=\max(0,\varphi(z))$, $\varphi_{-1}
:=d^{-1}\varphi_0 \circ f_0$, $E_0:=\varphi_0^{-1}(0)$
et $E_{-1}:=f_0^{-1}(E_0)$. Alors $\varphi_0$
est subharmonique; $\varphi_{-1}={d'}^{-1}\varphi_0\circ g_0$ et
$E_{-1}=g_0^{-1}(E_0)$. Comme $\varphi$ tend vers l'infini quand
$z\rightarrow \infty$, $E_0$ est compact.
Alors $\varphi_0$ est la
fonction de Green de $\P^1\setminus E_0$ avec un seul p\^ole en
$\infty$. On en d\'eduit que $E_0$ est de capacit\'e
logarithmique positive ({\it voir} par exemple, \cite[III.8]{Tsuji}). 
\par
Notons $\Sigma(d,d',\alpha)$
l'ensemble des
couples $(f,g)$ o\`u $f$ (resp. $g$) est un polyn\^ome de degr\'e
$d$ (resp. $d'$) \`a coefficient dominant $1$ (resp. $\alpha$)
qui s'annule en $0$.
\begin{lemme} Il existe un couple unique
$(f_1,g_1)\in\Sigma(d,d',\alpha^d)$
et un compact $E_1$ de capacit\'e logarithmique positive
tels que $f_1\circ g_0=g_1\circ f_0$ et tels que
$E_0=f_1^{-1}(E_1)=g_1^{-1}(E_1)$. De plus, on a ${\cal
  C}_{f_1}=g_0({\cal C}_{f_0})$ et ${\cal C}_{g_1}=f_0({\cal
  C}_{g_0})$.  
\end{lemme}
\begin{preuve} D'apr\` es le corollaire 3 et 
la proposition 2, il existe un polyn\^ome
  $\Phi$ de degr\'e $dd'$ et des polyn\^omes $f_1$ et $f_2$ tels que 
$\Phi=f_1\circ g_0=g_1\circ f_0$. Quitte \`a remplacer $\Phi$, $f_1$
et $g_1$  par
$\sigma\circ \Phi$, $\sigma\circ f_1$ et $\sigma\circ g_1$, on peut
supposer que $\Phi(0)=0$ et que le coefficient dominant de $\Phi$ soit
$\alpha^d$ o\`u $\sigma$ est un certain polyn\^ome lin\'eaire. 
On a alors $(f_1,g_1)\in\Sigma(d,d',\alpha^d)$.
\par
Montrons qu'au voisinage de $\infty$,  
$\delta_{f_1}$ pr\'eserve les lignes de niveau de
$\varphi_0$. Soient $a_1$ et $a_2$ suffisamment proches de
$\infty$ tels que $f_1(a_1)=f_1(a_2)$. Il faut prouver que
$\varphi_0(a_1)=\varphi_0(a_2)$. Il existe $b_1$ et $b_2$ tels
que $g_0(b_1)=a_1$ et $g_0(b_2)=(a_2)$. Alors
$\Phi(b_1)=\Phi(b_2)$. Par construction de $\Phi$ ({\it voir} la
preuve de la proposition 2), il existe $m$
et $n$ tels que $b_1=\delta_{f_0}^m\circ\delta_{g_0}^n(b_2)$. Comme
$\varphi_{-1}=d^{-1}\varphi_0\circ f_0 ={d'}^{-1}\varphi_0\circ
g_0$, les applications $\delta_{f_0}$ et $\delta_{g_0}$
pr\'eservent les
lignes de niveau de $\varphi_{-1}$. D'o\`u
$\varphi_{-1}(b_1)=\varphi_{-1}(b_2)$. On obtient
$$\varphi_0(a_1)=\varphi_0\circ g_0(b_1)=d'\varphi_{-1}(b_1)
=d'\varphi_{-1}(b_2)=\varphi_0\circ g_0(b_2)=\varphi_0(a_1).$$
Alors au voisinage de $\infty$, $\delta_{f_1}$ pr\'eserve les
lignes de niveau de $\varphi_0$, i.e. les lignes de niveau de
$\varphi_0$ sont r\'eunions de fibres de $f_1$. Comme $\varphi_0$ est
harmonique dans $\C\setminus E_0=\C\setminus \varphi_0^{-1}(0)$, elle
est r\'eelle analytique dans $\C\setminus E_0$. Par analyticit\'e,
$\C\setminus E_0$ est une r\'eunion de fibres de $f_1$. 
Par cons\'equent, $E_0$ est une
r\'eunion de fibres de $f_1$. Posons
$E_1:=f_1(E_0)$. Alors $E_0=f_1^{-1}(E_1)$.
Il est clair que $E_1$ est de
capacit\'e logarithmique positive. Les relations $f_1\circ
g_0=g_1\circ f_0$ et $f_0^{-1}(E_0)=g_0^{-1}(E_0)$
entra\^{\i}nent $g_1^{-1}(E_1)=E_0$.  Les polyn\^omes $f_1$ et
$g_1$ sont uniques car la fonction $\Phi$ est unique ({\it voir} la
preuve de la proposition 2). D'apr\`es la proposition 2, on a ${\cal
  C}_{f_1}=g_0({\cal C}_{f_0})$ et ${\cal C}_{g_1}=f_0({\cal
  C}_{g_0})$.   
\end{preuve}
\begin{remarque} {\rm
1. On peut construire les couples
$(f_k,g_k)\in\Sigma(d,d',\alpha^{d^k})$ et les compacts $E_k$
tels que $f_k\circ
g_{k-1}=g_k\circ f_{k-1}$ et
$E_{k-1}=f_k^{-1}(E_k)=g_k^{-1}(E_k)$.
\par 
2. On fixe un $k\geq 1$ et un
$m\geq 0$. On
pose $\tilde d:=d^k$, $\tilde
d':={d'}^k$, $\tilde
\alpha:=\alpha^{d^m(d^{k-1}+d^{k-2}d'+\cdots +{d'}^{k-1})}$,
$\tilde f_i:=f_{ik+k+m-1}\circ\cdots\circ f_{ik+m}$, $\tilde
g_i:=g_{ik+k+m-1}\circ\cdots \circ g_{ik+m}$,
$\tilde E_i:=E_{ik+k+m-1}$ pour $i=0$ ou 1. Alors 
$(\tilde f_0, \tilde
g_0)\in \Sigma(\tilde d,\tilde d',\tilde\alpha)$ et
$(\tilde f_1, \tilde
g_1)\in \Sigma(\tilde d,\tilde d',\tilde\alpha^{\tilde d})$. On
v\'erifie facilement que $\tilde f_1\circ
\tilde g_0 =\tilde g_1\circ \tilde f_0$ et que
$\tilde E_0=\tilde f_1^{-1}(\tilde E_1)=\tilde
g_1^{-1}(\tilde E_1)$. Par l'unicit\'e, $(\tilde
f_1, \tilde g_1)$ est le couple que l'on peut construire
comme dans le lemme 1 mais pour les polyn\^omes $\tilde f_0$ et
$\tilde g_0$.
}
\end{remarque}
On remarque qu'un couple $(f,g)\in \Sigma(d,d',\alpha)$ est
d\'etermin\'e uniquement par les points critiques de $f$ et
de $g$. Notons $\Pi_{d,d',\alpha}:
\C^{d-1}\times \C^{d'-1}\longrightarrow \Sigma(d,d',\alpha)$
l'application qui associe
un point $(x,y)=(x_1,\ldots,x_{d-1},y_1,\ldots, y_{d'-1})$ le couple
$(f,g)\in\Sigma(d,d',\alpha)$ v\'erifiant ${\cal
C}_f=\{x_1,\ldots,x_{d-1}\}$ et ${\cal C}_g=\{y_1,\ldots,
y_{d'-1}\}$. Cette application d\'efinit
un rev\^etement ramifi\'e au-dessus de $\Sigma(d,d',\alpha)$.
On d\'efinit
l'application
$\D_{d,d'}:\C^{d-1}\times\C^{d'-1}\times\C^* \longrightarrow
\C^{d-1}\times\C^{d'-1}\times \C^*$ par:
$$\D_{d,d'}(x,y,\alpha):
=(g(x_1),\ldots, g(x_{d-1}),f(y_1),\ldots,f(y_{d'-1}),\alpha^d)$$
o\`u $(f,g):=\Pi_{d,d',\alpha}(x,y)$. Les polyn\^omes $f$ et $g$
sont d\'etermin\'es par les formules explicites suivantes:
$$f(z)=\frac{1}{d}\int_0^z (t-x_1)\ldots (t-x_{d-1})dt$$ et
$$g(z)=\frac{\alpha^d}{d'}\int_0^z (t-y_1)\ldots (t-y_{d'-1})dt.$$
Il est clair que $\D_{d,d'}$ est un endomorphisme
polynomial. 
\begin{remarque} \rm
 1. D'apr\`es le lemme
pr\'ec\'edent, si $\Pi_{d,d',\alpha}(x,y)=(f_0,g_0)$ on a
$\Pi_{d,d',\alpha^d}(x^*,y^*)=(f_1,g_1)$ o\`u
$(x^*,y^*,\alpha^d):=\D_{d,d'}(x,y,\alpha)$.
\par
2. D'apr\`es la remarque 1, si $\Pi_{\tilde d,\tilde d',\tilde
\alpha}(\tilde x, \tilde y)=(\tilde f_0,\tilde g_0)$ on a
$\Pi_{\tilde d,\tilde d',\tilde \alpha}(\tilde x^*,\tilde
y^*)=(\tilde f_1,\tilde g_1)$ o\`u $(\tilde x^*,\tilde
y^*,\tilde\alpha^{\tilde d}):=\D_{\tilde d,\tilde d'}(\tilde
x,\tilde y,\tilde \alpha)$.
\end{remarque}
\section{Ensemble invariant ${\cal N}(d,d')$}
Notons ${\cal M}(d,d')$ l'ensemble des points
$(x,y,\alpha)\in\C^{d-1}\times \C^{d'-1}\times\C^*$ v\'erifiant
$f^*\circ g=g^*\circ f$ o\`u 
$(f,g):=\Pi_{d,d',\alpha}(x,y)$,
$(x^*,y^*,\alpha^d):=\D_{d,d'}(x,y,\alpha)$ et
$(f^*,g^*):=\Pi_{d,d',\alpha^d}(x^*,y^*)$. Notons ${\cal N}(d,d')$
l'ensemble des $(x,y,\alpha)\in{\cal M}(d,d')$ v\'erifiant les
deux propri\'et\'es suivantes:
\begin{trivlist}
\item[]1. ${\cal P}(d,d')$: pour tout $n\geq 0$, on a
$\D^n_{d,d'}(x,y,\alpha)\in {\cal M}(d,d')$.
\item[]2. Pour tous $k\geq 1$ et $m\geq 0$, si $\Pi_{\tilde
d,\tilde d',\tilde \alpha}(\tilde x,\tilde y)=
(\tilde f,\tilde g)$ alors $(\tilde x,\tilde y,\tilde \alpha)$
v\'erifie la condition ${\cal P}(\tilde d,\tilde d')$ o\`u
$(x_n,y_n,\alpha_n):=\D_{d,d'}^n(x,y,\alpha)$,
$(f_n,g_n):=\Pi_{d,d',\alpha_n}(x_n,y_n)$, $\tilde
f:=f_{k+m-1}\circ\cdots\circ f_m$, $\tilde
g:=g_{k+m-1}\circ\cdots\circ g_m$, $\tilde d:=d^k$, $\tilde
d':={d'}^k$ et
$\tilde\alpha:=\alpha^{d^m(d^{k-1}+d^{k-2}d'+\cdots+
{d'}^{k-1})}$ ({\it voir} la remarque 2).  
\end{trivlist}
Alors ${\cal N}(d,d')$ est un sous-ensemble
alg\'ebrique {\it faiblement invariant} par
$\D_{d,d'}$ i.e. $\D_{d,d'}({\cal N}(d,d'))\subset {\cal
N}(d,d')$. De plus, $\D_{d,d'}^n({\cal N}(d,d'))$ est faiblement
invariant par $\D_{d,d'}$ pour tout $n\geq 0$.
\par
Soient $f_0$, $g_0$, $\alpha$ et $E_0$ v\'erifiant les
hypoth\`eses du paragraphe pr\'ec\'edent. D'apr\`es le lemme
1 et les remarques 1, 2, on a $(x,y,\alpha)\in {\cal
N}(d,d')$ pour tout $(x,y)$ v\'erifiant
$\Pi_{d,d',\alpha}(x,y)=(f_0,g_0)$.
\par
Nous construisons maintenant deux sous-ensembles
${\cal C}_1(d,d')$ et ${\cal C}_2(d,d')$
de ${\cal N}(d,d')$ gr\^ace \`a des exemples pr\'ecis sur
$f_0$, $g_0$, $\alpha$ et $E_0$. Par suite, on montre que
${\cal N}(d,d')={\cal C}_1(d,d')\cup{\cal C}_2(d,d')$.
\par
Soient $\sigma_1$, $\sigma_2$ deux
polyn\^omes lin\'eaires, $a\not=0$ et $\alpha\not=0$ tels que
$(f,g)\in \Sigma(d,d',\alpha)$ o\`u  
$f(z):=\sigma_1\circ (z^d)\circ\sigma_2$ et 
$g(z):=\sigma_1\circ (az^{d'})\circ\sigma_2$. On a
$f^{-1}(E)=g^{-1}(E)$ pour $E:=\sigma_1(\{z\in\C: |z|=|a|^{d/(d-d')}\})$.
Par cons\'equent, 
$(x,y,\alpha)\in{\cal N}(d,d')$ pour tout $(x,y)$ v\'erifiant  
$\Pi_{d,d',\alpha}(x,y)=(f,g)$. On note
${\cal C}_1(d,d')$ l'ensemble de ces points $(x,y,\alpha)$. 
\par
Notons $\T_k$ le polyn\^ome de Tchebychev de degr\'e $k$
d\'efini par $\T(\cos z):=\cos kz$. On sait que l'ensemble de
Julia de $\T_k$ est l'intervalle $[-1,1]$, que le coefficient
dominant de $\T_k$ est \'egal \`a $2^{k-1}$ et que les points
critiques de $\T_k$ sont les points $\cos t\not=\pm 1$ avec $t\in\R$
v\'erifiant $\sin kt=0$. 
Soient $\sigma_1$, $\sigma_2$ deux polyn\^omes lin\'eaires et
$\alpha\in\C^*$ tels
que $(f,g)\in\Sigma(d,d',\alpha)$ o\`u
$f:=\sigma_1\circ(\pm \T_d)\circ\sigma_2$,
$g:=\sigma_1\circ(\pm\T_{d'})\circ\sigma_2$. Posons 
$E:=\sigma_1([-1,1])$. Alors
$f^{-1}(E)=g^{-1}(E)$.
On en d\'eduit que $(x,y,\alpha)\in{\cal N}(d,d')$.
Notons ${\cal C}_2(d,d')$
l'ensemble de tels points $(x,y,\alpha)$.
\begin{lemme} ${\cal C}_1(d,d')$ et ${\cal C}_2(d,d')$ sont des
courbes alg\'ebriques r\'eductibles dont aucune composante n'est
incluse dans un hyperplan du type $\{\alpha=\mbox{constante}\}$.
\end{lemme}
\begin{preuve} 
Soient
$\sigma_1(z)=a_1z+b_1$, $\sigma_2(z)=a_2z+b_2$ et $f$, $g$, $\alpha$,
$x$, $y$
d\'efinis ci-dessus.
\par
Pour la courbe ${\cal C}_1(d,d')$,
comme $(f,g)\in\Sigma(d,d',\alpha)$, on
obtient les relations suivantes: $a_1a_2^d=1$,
$a_1aa_2^{d'}=\alpha$ et $a_1b_2^d+b_1=a_1ab_2^{d'}+b_1=0$. On a
\'egalement, $x=(-b_2/a_2,\ldots,-b_2/a_2)$ et
$y=(-b_2/a_2,\ldots,-b_2/a_2)$. On obtient facilement que
$b_2=0$ ou $\alpha=(b_2/a_2)^{d-d'}$. Ceci montre que ${\cal
C}_1$ est une courbe alg\'ebrique 
r\'eductible et qu'aucune de ses composantes
n'est incluse dans $\{\alpha=\mbox{constante}\}$.
\par
Pour la courbe ${\cal C}_2(d,d')$, on obtient $a_12^{d-1}a_2^d=1$,
$a_12^{d'-1}a_2^{d'}=\alpha$ et
$\pm a_1 \T_d(b_2)+b_1=\pm a_1\T_{d'}(b_2)+b_1=0$. On a
$\{x_1,\ldots,x_{d-1}\}=\sigma_2^{-1}({\cal C}_{\T_{d_2}})$,
  $\{y_1,\ldots,y_{d'-1}\}=\sigma_2^{-1}({\cal C}_{\T_{d_2}})$. On remarque
  que $b_2$ est une solution de l'\'equation $\pm
  \T_d(z)=\pm\T_{d'}(z)$. Cette \'equation n'a qu'un nombre fini de
  solution. Les autres nombres $a_1$ et $a_2$ (resp. $b_1$) s'\'ecrivent en
  fonction de $\alpha$ (resp. de $\alpha$ et de $b_2$). 
Par cons\'equent, ${\cal
    C}_2(d,d')$ est une courbe alg\'ebrique 
dont aucune composante n'est incluse
  dans un hyperplan du type $\{\alpha=\mbox{constante}\}$.  
\end{preuve}
\begin{lemme} 
\begin{trivlist}
\item[] {\rm 1.}
$\D_{d,d'}^{-1}({\cal C}_1(d,d'))\cap{\cal N}(d,d')\subset {\cal
C}_1(d,d')$. 
\item[] {\rm 2.}
$\D_{d,d'}^{-1}({\cal C}_2(d,d'))\cap{\cal N}(d,d')\subset {\cal
C}_2(d,d')$. 
\end{trivlist}  
\end{lemme}
\begin{preuve} {\bf 1.} Soit $(x,y,\alpha)\in
\D_{d,d'}^{-1}({\cal C}_1(d,d'))\cap{\cal N}(d,d')$. Posons
$(f,g):=\Pi_{d,d',\alpha}(x,y)$,
$(x_1,y_1,\alpha_1):=\D_{d,d'}(x,y,\alpha)$ et
$(f_1,g_1):=\Pi_{d,d',\alpha_1}(x_1,y_1)$. Comme
$(x_1,y_1,\alpha_1)\in{\cal C}_1(d,d')$, il existe des
polyn\^omes lin\'eaires $\sigma_1$, $\sigma_2$ et une constante
non nulle $a$ tels que $f_1(z)=\sigma_1([\sigma_2(z)]^d)$ et
$g_1=\sigma_1(a[\sigma_2(z)^{d'}])$. Comme $(x,y,\alpha)\in{\cal
N}(d,d')$, on a $f_1\circ g=g_1\circ f$. Posons
$f^*:=\sigma_2\circ f$ et
$g^*:=\sigma_2\circ g$. Alors
$a[f^*(z)]^{d'}=[g^*(z)]^d$.      
Posons $\Phi(z):=a[f^*(z)]^{d'}=[g^*(z)]^d$.
Soit $\lambda$ une racine de $\Phi$. Alors la multiplicit\'e de
$\lambda$ est divisible par $d$ et par $d'$. Comme $d$ et $d'$ sont
premiers entre eux, la multiplicit\'e de $\lambda$ est divisible par
$dd'$. D'autre part, $\deg \Phi=dd'$. On d\'eduit que $\lambda$
est la
seule racine de $\Phi$. Il est \'egalement la seule racine
de $f^*$ et de $g^*$. Alors il existe un polyn\^ome lin\'eaire 
$\sigma$ et un $b\in\C$ tels que
$f^*(z)=[\sigma(z)]^d$ et $g^*(z)=b[\sigma(z)]^{d'}$. D'o\`u 
$(x,y,\alpha)\in {\cal C}_1(d,d')$. 
\par
{\bf 2.} De m\^eme mani\`ere, on se ram\`ene \`a une
\'equation du type $\T_d\circ g^*=\pm \T_{d'}\circ f^*$. Il
faut montrer qu'il existe un polyn\^ome lin\'eaire $\sigma$ tel
que $f^*=\pm \T_d\circ \sigma$ et $g^*=\pm \T_{d'}\circ\sigma$.
Il est clair que
${f^*}^{-1}([-1,1])= {g^*}^{-1}([-1,1])$. On d\'eduit de la
d\'efinition de $\D_{d,d'}$
que $g^*({\cal C}_{f^*})={\cal C}_{\T_d}$. Par
cons\'equent, les points critiques de $f^*$ sont tous de
multiplicit\'e 1. De m\^eme pour $g^*$. 
\par 
Pour 
$|z|$ suffisamment grand on a $\F_{f^*}(z)\cap
\F_{g^*}(z)=\{z\}$. En effet, utilisant les 
d\'eveloppements asymptotiques de
$\delta_{f^*}$ et $\delta_{g^*}$, on obtient pour tous $1\leq
n\leq d-1$ et $1\leq m\leq d'-1$:
$$\lim_{z\rightarrow
  \infty}\frac{\delta_{f^*}^m(z)}{\delta_{g^*}^n(z)}
=\exp(2m\pi i/d-2n\pi i/d')\not =0$$
car $d$ et $d'$ sont premiers entre eux. 
Par analyticit\'e,  pour un $z$ g\'en\'erique $\F_{f^*}(z)\cap
\F_{g^*}(z)=\{z\}$.
Soit $p\in{\cal C}_{f^*}$. Montrons que $f^*(p)=\pm 1$. Supposons
que $f^*(p)=a\not=\pm 1$. On sait que $g^*({\cal C}_{f^*})={\cal
C}_{\T_d}\subset [-1,1]$ et
${f^*}^{-1}([-1,1])={g^*}^{-1}([-1,1])$. D'o\`u $a\in ]-1,1[$.
Alors au voisinage de $p$, ${f^*}^{-1}([-1,1])$ est la r\'eunion
de deux courbes r\'eelles analytiques qui se coupent en $p$. 
Par cons\'equent, $p$ est
un point critique de $g^*$. Comme $\delta_{f^*}$ et
$\delta_{g^*}$ commutent, leurs prolongements analytiques
commutent aussi au voisinage de $p$. On en d\'eduit que
$\F_{f^*}(z)\cap \F_{g^*}(z)$ contient au moins deux \'el\'ements pour
tout $z$ suffisamment proche de $p$. 
C'est une contradiction. Donc $f^*(p)=\pm 1$. De
m\^eme $g^*(q)=\pm 1$ pour $q\in{\cal C}_{g^*}$.
\par
Le fait que $f^*$ est de degr\'e $d$ implique 
que pour $d$ impair ${f^*}^{-1}(1)$
est une
r\'eunion de $(d-1)/2$ points critiques et d'un point non
critique; pour $d$ pair ${f^*}^{-1}(1)$ est la r\'eunion de
$d/2$ points critiques ou la r\'eunion
de $d/2-1$ points critiques avec deux
points non critiques. De m\^eme pour ${f^*}^{-1}(-1)$. Quitte \`a
remplacer $f^*$, $g^*$ par $\pm f^*\circ\sigma$ et $g^*\circ
\sigma$ pour un certain polyn\^ome lin\'eaire $\sigma$, on peut
supposer que $\pm 1$ ne sont pas critiques pour $f$ et que
pour $d$ impair $f^*(1)=1$, $f^*(-1)=-1$ et pour $d$
pair $f^*(1)=f^*(-1)=1$. On remarque qu'au voisinage de $\pm
1$, ${f^*}^{-1}([-1,1])$ est un arc r\'eel analytique. Par
cons\'equent, $g^*(\pm 1)=\pm 1$ et $\pm 1$ ne sont pas critiques
pour $g^*$. Quitte \`a remplacer $g^*$ par $\pm g^*$, on peut
supposer que pour $d'$ impair $g^*(1)=1$, $g^*(-1)=-1$ et pour $d'$
pair $g^*(1)=g^*(-1)=1$.
\par
Alors il existe des polyn\^omes $P$, $Q$ tels que pour $d$
impair $f^*(z)+1=(z+1)P^2(z)$, $f^*(z) -1=(z-1)Q^2(z)$
et pour $d$ pair
$f^*(z) +1=P^2(z)$, $f^*(z)-1=(z-1)(z+1)Q^2(z)$. Posons
$\psi(z):=(z+z^{-1})/2$. On v\'erifie facilement qu'il existe
une fonction rationnelle $R(z)$ telle que:
$$\frac{f^*+1}{f^*-1}\circ\psi(z)=R^2(z).$$
On en d\'eduit que $f^*(z)\circ \psi=(F+F^{-1})/2$ o\`u
$F:=(R+1)/(R-1)$. On 
a donc $\psi^{-1}\circ f^*\circ \psi=F^{\pm 1}$. Ceci
implique que $\deg F=d$. De m\^eme,
il existe une fonction rationnelle $G$ de degr\'e $d'$ telle que
$\psi^{-1}\circ g^*\circ \psi=G^{\pm 1}$. D'autre part,
$\psi^{-1}\circ \T_d\circ\psi(z)=z^{\pm d}$ et
$\psi^{-1}\circ\T_{d'}\circ \psi(z)=z^{\pm d'}$. On
d\'eduit de la relation $\T_d\circ g^*=\pm\T_{d'}\circ f^*$
que $F^{\pm d'}=\pm G^{\pm d}$. Alors comme dans la partie
pr\'ec\'edente, les multiplicit\'es
des z\'eros et des p\^oles de $F^{\pm d'}=\pm G^{\pm
d}$ sont divisibles par $dd'$. Or c'est une
fonction de degr\'e $dd'$. D'o\`u $F(z)=Az^{\pm d}$ et
$G(z)=Bz^{\pm d'}$. Le fait que $f^*(1)=g^*(1)=1$ entra\^{\i}ne
$A=B=1$. D'o\`u $f^*=\T_d$ et $g^*=\T_{d'}$.       
\end{preuve}
\begin{lemme} Soit $(x,y,\alpha)\in{\cal N}(d,d')$ un point
pr\'ep\'eriodique de $\D_{d,d'}$, i.e. ${\cal
D}_{d,d'}^{k+m}(x,y,\alpha)=\D_{d,d'}^m(x,y,\alpha)$ pour
certains $k\geq 1$ et $m\geq 0$. Alors 
$(x,y,\alpha)$ appartient \`a
${\cal C}_1(d,d')\cup{\cal C}_2(d,d')$.
\end{lemme}
\begin{preuve} 
On utilise les notations de la
d\'efinition de l'ensemble ${\cal N}(d,d')$.
Soient
$(\tilde x^*, \tilde y^*,\tilde \alpha^*):=\D_{\tilde d,\tilde
d'}(\tilde x,\tilde y,\tilde \alpha)$ et $(\tilde f^*,\tilde
g^*):=\Pi_{\tilde d,\tilde d',\tilde\alpha}(\tilde x^*,\tilde
y^*)$. Comme
$\D_{d,d'}^{k+m}(x,y,\alpha)=\D_{d,d'}^m(x,y,\alpha)$, on a
$\tilde f^*=\tilde f$ et $\tilde g^*=\tilde g$ ({\it voir} la remarque
2). Par d\'efinition
de ${\cal N}(\tilde d,\tilde d')$,
on a $\tilde f^*\circ\tilde g=\tilde g^*\circ
\tilde f$. D'o\`u $\tilde f\circ\tilde g=\tilde g\circ\tilde f$.
Cette \'equation a \'et\'e r\'esolue par Fatou et Julia
\cite{Fatou,Julia,Ritt,Eremenko}. Dans notre cas, $\tilde d>1$ et $\tilde
d'>1$ sont premiers entre eux. D'apr\`es le th\'eor\`eme de Fatou-Julia, 
il existe un polyn\^ome lin\'eaire
$\sigma_1$ tel que l'une des conditions suivantes soit vraie:
\begin{enumerate}
\item $\sigma_1\circ \tilde f\circ\sigma_1^{-1}= z^{\tilde d}$ et
$\sigma_1\circ \tilde g \circ\sigma_1^{-1}=az^{\tilde d'}$ o\`u
$a\not =0$ est
une constante.
\item $\sigma_1\circ \tilde f\circ\sigma_1^{-1}=
\pm \T_{\tilde d}$ et
$\sigma_1\circ \tilde g \circ\sigma_1^{-1}=\pm \T_{\tilde d'}$.
\end{enumerate}
Consid\'erons le second cas, le premier cas sera trait\'e de
m\^eme mani\`ere. On remarque que $\T_{rs}=\T_r\circ \T_s$. En
particulier, $\T_{\tilde d}=\T_{\tilde d/d}\circ\T_d$.  
Le fait que $\tilde f=f_{k+m-1}\circ \cdots \circ f_m$ implique 
$$(\sigma_1\circ 
f_{k+m-1}\circ \sigma_1^{-1})\circ \cdots (\sigma_1 \circ f_m\circ
\sigma_1^{-1})=\sigma_1\circ \tilde f\circ \sigma^{-1}=\pm \T_{\tilde d}.$$
D'apr\`es la proposition 1, il existe un polyn\^ome lin\'eaire
$\sigma_2$ tel que $\sigma_1\circ f_m\circ \sigma_1^{-1}=\sigma_2
\circ\T_d$. De
m\^eme, il existe $\sigma_2'$ tel que $\sigma_1\circ
g_m\circ\sigma_1^{-1}=\sigma_2'\circ \T_{d'}$. Alors 
$f_m=\sigma_3\circ
\T_d\circ\sigma_1$ et
$g_m=\sigma_3'\circ\T_{d'}\circ\sigma_1$ o\`u
$\sigma_3:=\sigma_1^{-1}\circ\sigma_2$ et
$\sigma_3':=\sigma_1^{-1}\circ\sigma_2$. On sait que pour $d$ et $d'$
premiers entre eux l'ensemble
critique de $\T_d$ (resp. de $\T_{d'}$) est invariant par $\T_{d'}$
(resp. par $\T_d$). Par construction de
$\D_{d,d'}$, on a $${\cal C}_{f_{m+1}}=g_m({\cal
  C}_{f_m})=\sigma_3'\circ \T_{d'}({\cal C}_{\T_d})=
\sigma_3'({\cal C}_{\T_d}).$$
Alors il existe un polyn\^ome lin\'eaire $\sigma_4$ tel que
$f_{m+1}=\sigma_4\circ\T_d\circ{\sigma_3'}^{-1}$. De m\^eme, il existe
un polyn\^ome lin\'eaire $\sigma_4'$ tel que
$g_{m+1}=\sigma_4'\circ \T_{d'}\circ\sigma_3^{-1}$.
\par
Comme on a montr\'e ci--dessus pour $f_m$ et $g_m$, il suffit de remplacer
$m$ par $m+1$ afin d'obtenir
$f_{m+1}=\sigma_6\circ\T_d\circ\sigma_5$ et
$g_{m+1}=\sigma_6'\circ\T_d\circ \sigma_5$ o\`u $\sigma_5$,
$\sigma_6$ et $\sigma_6'$ sont lin\'eaires. On d\'eduit des
quatres derni\`eres \'egalit\'es que l'ensemble critique de $\T_d$
(resp. de $\T_{d'}$) est
invariant par $\sigma_5\circ\sigma_3'$ (resp. par
$\sigma_5\circ\sigma_3$).
D'o\`u $\sigma_5\circ\sigma_3' (z)=\pm z$ et
$\sigma_5\circ\sigma_3(z)=\pm z$. Par cons\'equent,
$\sigma_3'(z)=\sigma_3(\pm z)$ et donc $f_m=\sigma_3\circ \T_d\circ
\sigma_1$ et $g_m=\sigma_3\circ (\pm
\T_{d'}) \circ \sigma_1$. Ceci signifie que
$(f_m,g_m,\alpha_m)\in{\cal C}_2(d,d')$. D'apr\`es le lemme 3,
$(f,g,\alpha)\in {\cal C}_2(d,d')$ car
$(f,g,\alpha)=(f_0,g_0,\alpha_0)$. 
\end{preuve}
\begin{proposition} On a 
${\cal N}(d,d')={\cal C}_1(d,d')\cup{\cal C}_2(d,d')$.
\end{proposition}
Soit $S$ un sous-ensemble alg\'ebrique p\'eriodique de
${\cal N}(d,d')$, i.e. $\D_{d,d'}^n(S)=S$ pour un certain $n\geq 1$. 
On montre que $S\subset {\cal
C}_1(d,d')\cup{\cal C}_2(d,d')$.
Soit $a$ une racine d'ordre $d^{m-1}-1$ de l'unit\'e.
On pose $K_a$ l'ensemble des points $(x,y,a)$. Alors $K_a$ est
p\'eriodique de p\'eriode $m$.
\begin{lemme} Pour tout $a$, l'ensemble $S\cap K_a$ est fini.
\end{lemme}
\begin{preuve}
Soit $V$ une composante
irr\'eductible, p\'eriodique de $S\cap
K_a$. Il faut montrer que $\dim V=0$. Supposons par l'absurde que $\dim
V\geq 1$. 
Pour simplifier les notations, on suppose par la suite que
$a=1$ et on pose $\D:=\D_{d,d'}$. 
Notons $s:=(x,y)=(x_1,\ldots,x_{d-1},y_1,\ldots,y_{d'-1})$
les coordonn\'ees de
$K_a\simeq\C^{d+d'-2}$ . Comme $\D$ est polynomial,
elle se prolonge en une application m\'eromorphe de
$\P^{d+d'-2}$ dans lui-m\^eme. Notons encore $\D$ ce prolongement.
Soit $L:=\P^{d+d'-2}\setminus K_a$ l'hyperplan
\`a l'infini muni des coordonn\'ees homog\`enes
$w:=[x_1:\cdots:x_{d-1}:y_1:\cdots:y_{d'-1}]$. On pose
$(f,g):=\Pi_{d,d',a}(s)$. Les formules explicites des polyn\^omes
$f$ et $g$ sont donn\'ees
dans le paragraphe pr\'ec\'edent.
On remarque que $f(y_i)$ (resp. $g(x_i)$)
est un polyn\^ome homog\`ene de degr\'e $d$ (resp. $d'$) en
variables $x$ et $y$. Par cons\'equent, il existe une constante $c>0$
telle que $|f(y_j)|\leq c\lambda^d$ et $|g(x_i)|\leq
c\lambda^{d'}$ o\`u
$$\lambda:=\max\left(\max_{1\leq \nu\leq d-1}|x_\nu|,
\max_{1\leq\nu\leq d'-1} |y_\nu|\right).$$
Comme $d>d'$, l'ensemble
d'ind\'etermination $I$ de $\D$ est \'egal \`a
$$I=\{w\in L:\ f(y_1)=\cdots=f(y_{d'-1})=0\}$$
et l'ensemble  $X:=\D(L\setminus I)$ v\'erifie
$$X\subset \{w\in L:\ x_1=\cdots=x_{d-1}=0\}.$$
Comme $\dim V\geq 1$, l'intersection $\overline V\cap
L\not=\emptyset$. Comme $V$ est pr\'ep\'eriodique, $\overline V\cap (I\cup
X)\not=\emptyset$. 
Soient $s^{(n)}=(x^{(n)},y^{(n)})\in
V$ tendant vers un point $w_0\in \overline V\cap(I\cup X)$
quand $n\rightarrow +\infty$. On pose
$(f_n,g_n):=\Pi_{d,d',a}(s^{(n)})$, $\overline s^{(n)}:=\D(s^{(n)})$ et
$(\overline f_n,\overline g_n):=
\Pi_{d,d',a}(\overline s^{(n)})$.
Par d\'efinition de $\D_{d,d'}$, on a ${\cal C}_{\overline
f_n}=g_n({\cal C}_{f_n})$ et ${\cal C}_{\overline g_n}=f_n({\cal
C}_{g_n})$. 
Par d\'efinition de
${\cal N}(d,d')$,
on a $\overline f_n\circ g_n=\overline g_n\circ f_n$. 
Soient
$$\lambda_n=\max\left(\max_{1\leq\nu\leq d-1}|x^{(n)}_\nu|,
\max_{1\leq\nu\leq d'-1}|y^{(n)}_\nu|\right).$$
Alors $|\overline x^{(n)}_\nu|\leq c\lambda_n^{d'}$ et
$|\overline y^{(n)}_\nu|\leq c\lambda_n^{d}$.  
On pose $\sigma_1(z):=\lambda z$,
$\sigma_2(z):= \lambda^{d-1}z$, $\sigma_3(z):=\lambda^{d'-1}z$ et
$\sigma_4(z):=\lambda^{(d-1)(d'-1)}z$. On pose \'egalement
$f_n^*:=\sigma_2^{-1}\circ f_n\circ\sigma_1$,
$g_n^*:=\sigma_3^{-1}\circ g_n\circ \sigma_1$, $\overline
f_n^*:=\sigma_4^{-1}\circ \overline f_n\circ \sigma_3$ et $\overline
g_n^*:=\sigma_4^{-1}\circ \overline g_1\circ \sigma_2$. Alors
$\overline 
f_n^*\circ g_n^*=\overline g_n^*\circ f_n^*$, $(f_n^*,g_n^*)
\in \Sigma(d,d',1)$ et $(\overline f_n^*,\overline
g_n^*)\in \Sigma(d,d',1)$. On a aussi ${\cal C}_{f_n^*}=
\sigma_1^{-1}({\cal
C}_{f_n})=\sigma_1^{-1}\{x^{(n)}_1,\ldots,x^{(n)}_{d-1}\}$,
${\cal
C}_{g_n^*}=\sigma^{-1}\{y^{(n)}_1,\ldots,y^{(n)}_{d'-1}\}$,
${\cal C}_{\overline f^*}=\sigma_3^{-1}\{\overline
x^{(n)}_1,\ldots, \overline x^{(n)}_{d-1}\}$ et ${\cal
C}_{\overline g^*}= \sigma_2^{-1}\{\overline y^{(n)}_1,\ldots,
\overline y^{(n)}_{d'-1}\}$. Par d\'efinition de $\lambda_n$ et
des $\sigma_i$, les points critiques de $f^*$ et $g^*$ (resp. de
$\overline f^*$ et $\overline g^*$) sont de
modules major\'es par $1$ (resp. par $c$). De plus, au moins l'un des
points critiques de $f^*$ ou de $g^*$ est de module $1$. Le fait que
$(f^*,g^*)\in\Sigma(d,d',1)$ et $(\overline f^*,\overline g^*)\in
\Sigma(d,d',1)$ entra\^{\i}ne que les coefficients des
polyn\^omes $f_n^*$, $g_n^*$, $\overline f_n^*$ et $\overline g_n^*$ 
sont born\'es. On v\'erifie facilement que
${\cal C}_{\overline
f_n^*}=g_n^*({\cal C}_{f_n^*})$ et ${\cal C}_{\overline g_n^*}=
f_n^*({\cal C}_{g_n^*})$. 
Soient $F$, $G$,
$\overline F$, $\overline G$ quatre polyn\^omes tels que
$(F,G,\overline F,\overline G)$ soit adh\'erent \`a la suite
$(f_n^*,g_n^*,\overline f_n^*, \overline g_n^*)$. Par continuit\'e,
on a $\overline F\circ G=\overline G\circ F$,
${\cal C}_{\overline F}=
G({\cal C}_F)$ et ${\cal C}_{\overline G}=F({\cal C}_G)$. De plus, au
moins un point critique de $F$ ou de $G$ est de module 1.  
\par
{\it Cas 1.}--
Supposons que $w_0\in I$. On d\'eduit de la description de $I$
que $\lambda_n^{-d+1}\overline y^{(n)}_\nu$ tend vers
$0$ quand $n\rightarrow +\infty$. Ceci implique que les points
critiques de $\overline G$ sont tous nuls. D'o\`u $\overline
G(z)=z^{d'}$ et $[F(z)]^{d'}=\overline F\circ G(z)$ car $\overline
F\circ G=\overline G\circ F$.
Alors les
multiplicit\'es des z\'eros de $\overline F\circ G$ sont
divisibles par $d'$. Comme $d=\deg \overline F$ n'est pas
divisible par $d'$, il existe au moins une racine $a_1$ de $\overline F$
telle que sa multiplicit\'e $\alpha_1$ ne soit pas divisible par
$d'$.
\par
Supposons d'abord qu'il existe une autre racine $a_2$ de
$\overline F$ dont la multiplicit\'e $\alpha_2$ n'est pas
divisible par $d'$.
Soit $b_j$ un point arbitraire de 
$G^{-1}(a_j)$ \`a multiplicit\'e $\beta_j$.
Alors $\alpha_j\beta_j$ est divisible par $d'$.
Notons
$\alpha'_j$ le plus grand diviseur commun de $\alpha_j$ et $d'$.
Notons \'egalement $\nu_j=d'/\alpha'_j$.
Alors $\nu_j$ divise $\beta_j$.
On en d\'eduit qu'il
existe un polyn\^ome non constant
$K_j$ tel que $G(z)-a_j=[K_j(z)]^{\nu_j}$. 
Comme $\alpha_j$ ne divise pas $d'$,
on a $\nu_j\geq 2$. On obtient donc
$[K_1(z)]^{\nu_1}=[K_2(z)]^{\nu_2}-b^{\nu_2}$ o\`u
$b^{\nu_2}=a_2-a_1\not=0$. Ceci implique
$$[K_1(z)]^{\nu_1}=\prod_{j=0}^{\nu_2-1}[K_2(z)-\theta_jb].$$
o\`u $\theta_j:=\exp(2j\pi i/\nu_2)$.
Les facteurs du membre \`a droite sont deux \`a deux premiers entre
eux. 
Par cons\'equent, il existe des polyn\^omes $P_j$
tels que $K_2(z)-\theta_jb=[P_j(z)]^{\nu_1}$. On a
$$[P_1(z)]^{\nu_1}-[P_0(z)]^{\nu_1}=(\theta_1-\theta_0)b\not=0.$$
C'est une contradiction car le membre \`a gauche se factorise en
$\nu_1$ facteurs qui ne sont pas tous constants.
\par
Il reste le cas o\`u $a_1$ est la seule racine de $\overline F$
dont la multiplicit\'e n'est pas divivible par $d'$. Comme
$d=\deg \overline F$ et $d'$ sont premiers entre eux, $\alpha_1$
et $d'$ sont premiers entre eux. Par cons\'equent, tout point de
$G^{-1}(a_1)$ est de multiplicit\'e divisible par $d'$. Comme
$\deg G=d'$ et comme $G\in\Sigma(d,d',1)$, on a $G(z)=z^{d'}$. On
en d\'eduit que $a_1=0$. On peut donc \'ecrire $\overline
F(z)=z^{\alpha_1} [P(z)]^{d'}$ o\`u $P$ est un polyn\^ ome
unitaire. 
L'\'equation $[F(z)]^{d'}=\overline F\circ G(z)$ entra\^{\i}ne
$F(z)=z^{\alpha_1}P(z^{d'})$. Les \'egalit\'es suivantes sont obtenues
par les calculs de d\'eriv\'ees:
$$F'(z)=z^{\alpha_1-1}[\alpha_1 P(z^{d'})+d'z^{d'}P'(z^{d'})]$$
et
$$\overline F'(z)=z^{\alpha_1-1}[P(z)]^{d'-1}
[\alpha_1 P(z) +d'zP'(z)].$$
Soit $a$ une racine non nulle de $F'$, i.e. une racine de
$\alpha_1P(z^{d'})+d'z^{d'}P'(z^{d'})$. Alors $\exp(2k\pi i/d')a$ est
\'egalement une racine de $F'$ pour tout $0\leq k\leq d'-1$.
Comme ${\cal
C}_{\overline F}=G({\cal C}_F)$ et comme $G(z)=z^{d'}$,
toute racine de $\overline F'$ est du type $a^{d'}$, i.e. une
racine de $\alpha_1 P(z)+d'zP'(z)$. De plus, 
la multiplicit\'e de cette racine est divisible par $d'$. Soit
$b$ une racine de multiplicit\'e $nd'+m$ de $P(z)$ avec $0\leq
m\leq d'-1$. Alors $b$ est une racine de multiplicit\'e 
$nd'+m-1$ de $P'(z)$ et donc de
$\alpha_1 P(z)+d'zP'(z)$. Par cons\'equent,
$b$ est une racine de multiplicit\'e $(n'd'+m)(d'-1)+(nd'+m-1)$ de
$\overline F'$. Cette multiplicit\'e 
n'est pas divisible par $d'$. C'est impossible. D'o\`u $P(z)=1$ et
$F(z)=z^d$. C'est aussi une contradiction car au moins l'un des points
critiques de $F$ ou de $G$ est  de module 1. 
\par
{\it Cas 2.}-- Supposons maintenant que $w_0\in X$. Par la
description de $X$, $\lambda_n^{-1}x^{(n)}_\nu$ tend vers $0$
quand $n\rightarrow +\infty$. Par cons\'equent, les points
critiques de $f_n^*$ tendent vers $0$. On en d\'eduit que
$F(z)=z^d$ et que ${\cal C}_{\overline F}=G({\cal C}_F)=\{0\}$.
On a donc $\overline F(z)=z^{d}$.
On obtient alors $\overline G(z^d)=[G(z)]^d$. Ceci montre que les
racines de $\overline G(z^d)$ sont toutes de multiplicit\'e
divisible par $d$. En particulier, toute racine non nulle de
$\overline G$ est de multiplicit\'e divisible par $d$. Mais $\deg
\overline G=d'<d$. Donc $\overline G$ n'a pas de racine non nulle.
Alors $\overline G(z)=z^{d'}$ et donc $G(z)=z^{d'}$. C'est une
contradiction car au moins un point critique de $F$ ou de $G$ est de
module 1.   
\end{preuve}
{\it Fin de la preuve de la proposition 2.}--
Si $\dim S=0$, alors $S$ est simplement un point p\'eriodique.
D'apr\`es le lemme 4, $S\subset{\cal C}_1(d,d')\cup{\cal
C}_2(d,d')$. 
\par
Si $\dim S\geq 1$, d'apr\`es le lemme pr\'ec\'edent,
$S\cap K_a$ est un ensemble fini
pour tout $a$. Comme $S$ et $K_a$ sont p\'eriodiques, $S\cap K_a$
est p\'eriodique. Par cons\'equent, tout point de $S\cap K_a$ est
pr\'ep\'eriodique. D'apr\`es le lemme 4,   $S\cap
K_a\subset {\cal C}_1(d,d')\cup{\cal C}_2(d,d')$.
Comme $S$ est p\'eriodique et comme $\D_{d,d'}$ envoie
l'hyperplan $\{\alpha=c\}$ dans l'hyperplan $\{\alpha=c^d\}$,
$S\cap K_a$ est non vide sauf peut-\^etre pour un nombre fini de
$a$. On d\'eduit que $\dim S=1$ et que 
$S$ coupe ${\cal C}_1(d,d')\cup {\cal C}_2(d,d')$ en une
infinit\'e de points. D'o\`u $S\subset {\cal C}_1(d,d')\cup{\cal
C}_2(d,d')$. Ceci est vrai pour tout sous-ensemble alg\'ebrique
p\'eriodique de ${\cal N}(d,d')$. 
\par
Comme $\D_{d,d'}^n({\cal N}(d,d'))$
est faiblement invariant pour tout $n\geq 0$, toute composante de
${\cal N}(d,d')$ s'envoie par un $\D_{d,d'}^n$ dans une composante
p\'eriodique de ${\cal N}(d,d')$. Donc elle s'envoie par ${\cal
  D}_{d,d'}^n$ dans ${\cal C}_1(d,d')\cup{\cal C}_2(d,d')$. D'apr\`es
le lemme 3, elle est incluse dans ${\cal C}_1(d,d')\cup{\cal
  C}_2(d,d')$. 
\section{Preuves des th\'eor\`emes et remarques} 
{\it Preuve du th\'eor\`eme 1}--- Soient $\mu$, $f$ et $g$
v\'erifiant les hypoth\`eses du th\'eor\`eme 1. Si $d$ est
divisible par $d'$ ou si $d'$ est divisible par $d$, d'apr\`es le
corollaire 3, la condition 1 du th\'eor\`eme 1 est satisfaisante.
\par
Dans le cas contraire, d'apr\`es le corollaire 3, on peut
supposer que $d$ et $d'$ sont premiers entre eux et que $d>1$,
$d'>1$. Sans perdre en g\'en\'eralit\'e, on peut supposer que
$d>d'$. Alors d'apr\`es les paragraphes 3 et 4, il existe des
polyn\^omes lin\'eaires $\sigma_1$, $\sigma_2$ et un nombre
$\alpha\not=0$ tels que
$$(\sigma_1\circ
f\circ\sigma_2,\sigma_1\circ
g\circ\sigma_2)\in\Pi_{d,d',\alpha}({\cal N}(d,d')).$$
D'apr\`es la proposition 3, ${\cal N}(d,d')={\cal C}_1(d,d')\cup{\cal
C}_2(d,d')$. Par d\'efinition de ${\cal
C}_1(d,d')$ et ${\cal C}_2(d,d')$, il existe des polyn\^omes
lin\'eaires $\sigma_3$ et $\sigma_4$ tels que l'une des
conditions suivantes soit vraie:
\begin{enumerate}
\item $\sigma_3\circ f\circ\sigma_4(z)=z^d$ et $\sigma_3\circ
g\circ \sigma_4(z)=az^{d'}$ o\`u $a\not =0$ est une constante.
\item $\sigma_3\circ f\circ\sigma_4=\pm\T_d$ et $\sigma_3\circ
g\circ \sigma_4=\pm\T_{d'}$. 
\end{enumerate}
Posons $Q:=\sigma_3^{-1}\circ\sigma_4^{-1}$, $f_0:=f\circ Q^{-1}$ et
$g_0:=g\circ Q^{-1}$. On a $f=f_0\circ Q$ et $g=g_0\circ Q$. On a
aussi 
$f_0=\sigma_3^{-1}\circ z^d\circ \sigma_3$,
$g_0=\sigma_3^{-1}\circ (az^{d'})\circ \sigma_3$ ou
$f_0=\sigma_3^{-1}\circ (\pm\T_d)\circ \sigma_3$,
$g_0=\sigma_3^{-1}\circ (\pm\T_{d'})\circ \sigma_3$. Alors pour la
nouvelle coordonn\'ee $z':=\sigma_3^{-1}(z)$,
$f_0$ et $g_0$ v\'erifient la condition 2 ou la condition 
3 du th\'eor\`eme 1.
\par
\hfill $\square$ 
\begin{lemme} Soient $E\subset \C$ un compact, $f$ et $g$ deux
polyn\^omes de degr\'es $d>1$ et $d'>1$. Supposons que
$f^{-1}(E)=g^{-1}(E)$ et que $d$, $d'$ sont premiers entre eux.
\begin{enumerate}
\item Si $f(z)=az^d$ et $g(z)=bz^{d'}$ avec $a\not=0$ et
$b\not=0$, alors $E$ est
une r\'eunion de cercles centr\'es en $0$
\item Si $f(z)=\pm \T_d$ et $g(z)=\pm \T_{d'}$, alors $E=[-1,1]$.
\end{enumerate}
\end{lemme}
\begin{preuve} {\bf 1.} On note $S_r$ le cercle de centre $0$ et de rayon
  $r\geq 0$. Pour tout compact non vide 
  $K\subset \C$, on pose $A_K(r)$ le
  maximum des longueurs des composantes connexes de $S_r\setminus
  K$. On pose
$$A_K:=\sup\{A_K(r)\mbox{ pour tout } r>0 \mbox{ tel que } K\cap
S_r\not=\emptyset\}.$$
Posons $F:=f^{-1}(E)$. On a $A_F=d^{-1}A_E$. D'autre part,
$F=g^{-1}(E)$. D'o\`u $A_F={d'}^{-1}A_E$. On en d\'eduit que
$A_E=A_F=0$. Par cons\'equent, $E$ est une r\'eunion de cercles
centr\'es en $0$.   
\par
{\bf 2.} Notons $\varphi$ la fonction de Green de $\P^1\setminus
[-1,1]$ avec un seul p\^ole en $\infty$. On a
$d^{-1}\varphi\circ f={d'}^{-1}\varphi\circ g=\varphi$.
Notons $E_{-1}=f^{-1}(E)=g^{-1}(E)$. On a
$$\max_{E_{-1}}
\varphi(z)=d^{-1}\max_E\varphi(z)={d'}^{-1}\max_E\varphi(z).$$
Par cons\'equent, $\varphi(z)=0$ pour $z\in E_{-1}$. D'o\`u
$E_{-1}\subset [-1,1]$ et $E\subset [-1,1]$. Notons
$\psi(z):=(z+z^{-1})/2$. On a  $\psi^{-1}\circ f\circ
\psi(z)=\pm z^{\pm d}$, $\psi^{-1}\circ g\circ \psi(z)=\pm
z^{\pm d'}$. Posons $\tilde E:=\psi^{-1}(E)$. Alors
$\tilde E\subset\psi^{-1}([-1,1])=\{z: |z|=1\}$. On a $\tilde
f^{-1}(\tilde E)=\tilde g^{-1}(\tilde E)$ o\`u $\tilde f(z):=\pm
z^d$ et $\tilde g(z):=\pm z^{d'}$. D'apr\`es la partie
pr\'ec\'edente, $\tilde E$ est le cercle unit\'e. D'o\`u
$E=[-1,1]$.     
\end{preuve}
{\it Preuve du corollaire 1}--- Dans le corollaire 1, la
condition n\'ecessaire est \'evidente. Pour la condition
suffisante, supposons par l'absurde qu'il existe deux polyn\^omes
distints $f$ et $g$  tels que $f^{-1}(E)=g^{-1}(E)$.  
Notons $\Omega$ la composante connexe de
$\P^{-1}\setminus E$ qui contient $\infty$.
Alors $f^{-1}(\Omega)=g^{-1}(\Omega)$. Comme $E$
est de capacit\'e logarithmique positive, il existe une fonction
de Green $\varphi$ de $\Omega$
avec un seul p\^ole en $\infty$
\cite[III.8]{Tsuji}. Alors $d^{-1}\varphi\circ f$,
${d'}^{-1}\varphi\circ g$ sont les fonctions de Green de
$f^{-1}(\Omega)=g^{-1}(\Omega)$ avec un seul p\^ole en $\infty$.
Comme la fonction de Green est unique, on a $d^{-1}\varphi\circ
f={d'}^{-1}\varphi\circ g$. On pose $\varphi_0(z)=0$ si
$z\not\in\Omega$ et $\varphi_0(z)=\varphi(z)$ si $z\in\Omega$.
C'est une fonction subharmonique et $\mu:=i\partial\overline\partial
\varphi_0$ est la mesure d'\'equilibre de $E$ \cite[III]{Tsuji}. 
On obtient par les
relations pr\'ec\'edentes que $d^{-1}f^*(\mu)={d'}^{-1}g^*(\mu)$.   
D'apr\`es le th\'eor\`eme 1, on a $f=f_0\circ Q$ et $g=g_0\circ
Q$. D'o\`u $f_0^{-1}(E)=g_0^{-1}(E)$.
Si la condition 1 du th\'eor\`eme 1
est vraie, on a $f_0^{-1}(E)=g_0^{-1}(E)=E$. 
\par
Si la condition 2 du th\'eor\`eme 1 est vraie, 
d'apr\`es le lemme pr\'ec\'edent, $E$ est
une r\'eunion de cercles centr\'es en $0$. 
On a $P^{-1}(E)=E$ pour toute rotation $P$ de centre $0$.
\par
Si la condition 3 est vraie, d'apr\`es le lemme pr\'ec\'edent,
$E=[-1,1]$. Par cons\'equent, $P^{-1}(E)=E$
pour $P:=\T_k$.
\par
Dans les trois cas, on obtient une contradiction avec 
l'hypoth\`ese du corollaire 1.
\par\hfill $\square$
\\
{\it Preuve du corollaire 2.}--- Comme $E$ est de capacit\'e
logarithmique positive, $E$ est un ensemble infini. D'apr\`es le
corollaire 1, il existe un polyn\^ome $P\not=\id$
tel que $P^{-1}(E)=E$.
Si $P(z)=az+b$, on a $|a|=1$ et $a\not=1$ car $E$ est compact.
Alors $P$ est une rotation de centre $b/(1-a)$. C'est impossible.
On a donc $\deg P\geq 2$. On sait que $J_P$ est le plus petit
compact totalement invariant par $P$ qui contient plus qu'un
\'el\'ement. D'o\`u $J_P\subset E$. Comme $K_P$ est le
plus grand compact totalement invariant par $P$, on a $E\subset
K_P$.
\par
\hfill $\square$
\\
Dans le cas g\'en\'eral, si $E\subset \C$ est un compact et si $f$, $g$
sont deux polyn\^omes v\'erifiant $f^{-1}(E)=g^{-1}(E)$, il n'existe
pas de mesure $\mu$ \`a support dans $E$ telle que
$d^{-1}f^*(\mu)={d'}^{-1}g^*(\mu)$. Par exemple pour $E=\{0\}$,
$f(z)=z(z-1)$ et $g(z)=z^2(z-1)$, la seule mesure de
probabilit\'e $\mu$
support\'ee par $E$ est la masse de Dirac en $0$. On a
$f^{-1}(E)=g^{-1}(E)$ mais
$d^{-1}f^*(\mu)\not ={d'}^{-1}g^*(\mu)$. 
\begin{proposition} Soient $E\subset \C$ un compact, $f$ et $g$
deux polyn\^omes tels que $f^{-1}(E)=g^{-1}(E)$. Alors il existe
deux mesures de probabilit\'e $\mu_1$ et $\mu_2$ \`a support
dans $E$ telles que $g_*(d^{-1}f^*(\mu_1))=\mu_1$ et
$f_*({d'}^{-1}g^*(\mu_2))=\mu_2$.
\end{proposition}
\begin{preuve} Soit $\delta_0$
une mesure de probabilit\'e \`a support
dans $E$. Posons $\delta_n:=g_*(d^{-1}f^*(\delta_{n-1}))$ pour tout
$n\geq 1$. Alors $\delta_n$ est une mesure de probabilit\'e \`a
support dans $E$. Posons
$S_n:=(\delta_0+\cdots+\delta_{n-1})/n$. Alors $S_n$ est
\'egalement une mesure de probabilit\'e \`a support dans $E$.
Il existe une suite croissante $\{n_i\}_{i\in\N}$ telle que
$S_{n_i}$ tende faiblement vers une mesure $\mu_1$ quand
$i\rightarrow +\infty$ car l'ensemble des mesures de
probabilit\'e \`a support dans $E$ est compact. On en d\'eduit que
$g_*(d^{-1}f^*(S_{n_i}))$ tend faiblement vers
$g_*(d^{-1}f^*(\mu_1))$. D'autre part,
$g_*(d^{-1}f^*(S_{n_i}))-S_{n_i}=(\delta_{n_i}-\delta_0)/n_i$
tend vers
$0$. On obtient finalement $g_*(d^{-1}f^*(\mu_1))=\mu_1$. De
m\^eme pour $\mu_2$.  
\end{preuve}
\end{document}